%% file: main.tex
\documentclass[titlepage]{amsart}
\usepackage{amssymb}
\usepackage{amsthm}
\usepackage{bbm}
\usepackage{graphicx}

\usepackage[colorlinks=true]{hyperref}
\usepackage[all]{xypic}

\usepackage[final]{epsfig}
\usepackage{psfrag}
\usepackage{epstopdf}

\usepackage{amsaddr}

\newtheorem{thm}{Theorem}[section]

\newtheorem{prop}[thm]{Proposition}
\theoremstyle{definition}
\newtheorem{definition}{Definition}[section]

\theoremstyle{observation}

\theoremstyle{definition}

\newcommand{\voided}[1]{}

\newcommand{\defin}[1]{{\it #1}}
\newcommand{\R}{\mathbb{R}}
\newcommand{\N}{\mathbb{N}}
\newcommand{\Q}{\mathbb{Q}}

{\bf}{\it}
{\bf}{\it}
\newtheorem*{riemann-mapping-theorem}{Riemann Mapping  Theorem}{\bf}{\it}
{\bf}{\it}
{\bf}{\it}
{\bf}{\it}
{\bf}{\it}
{\bf}{\it}
{\bf}{\it}
{\bf}{\it}
{\bf}{\it}
{\bf}{\it}

\newenvironment{pf*}[1]{\proof[#1]}{\endproof}
\usepackage{euscript}

\newcommand{\cal}[1]{{\mathcal #1}}

\newcommand{\beq}{\begin{equation}}
\newcommand{\eeq}{\end{equation}}

\newtheorem{defn}{Definition}[section]

\newcommand{\tl}{\tilde}

\newcommand{\eps}{\epsilon}

\numberwithin{equation}{section}

\newcommand{\supp}{\operatorname{Supp}}

\newcommand{\cM}{{\mathcal M}}

\newcommand{\cP}{{\mathcal P}}

\renewcommand{\cD}{{\cal D}}

\newcommand{\RR}{{\mathbb R}}

\newcommand{\NN}{{\mathbb N}}

\newcommand{\QQ}{{\mathbb Q}}

\newcommand{\ignore}[1]{{}}

 \title[How to lose at Monte Carlo]{How to lose at Monte Carlo: a simple dynamical system whose typical statistical behavior is non computable}

\author{Cristobal Rojas}
\address{Departamento de Matem\'aticas, Universidad Andres Bello, \\Rep\'ublica 498, Santiago, Chile.\\ E-mail:  \texttt{\emph{crojas{@}mat-unab.cl}}}

\author{Michael Yampolsky}
\address{Department of Mathematics, University of Toronto,\\ 40 St George Street, Toronto, Ontario, Canada. \\ E-mail:  \texttt{\emph{yampol{@}math.toronto.edu}}}

 \thanks{ C.R was partially supported by the European Union's Horizon 2020 research and innovation program under the Marie Sklodowska-Curie grant agreement No 731143, by project FONDECYT  Regular  No 1190493 and project Basal PFB-03 CMM-Universidad de Chile.  M.Y. was partially supported by NSERC Discovery grant.}
\begin{document}

\begin{abstract}  We consider the simplest non-linear discrete dynamical systems, given by the logistic maps $f_{a}(x)=ax(1-x)$ of the interval $[0,1]$. We show that there exist real parameters  $a\in (0,4)$ for which almost every orbit of $f_a$ has the same statistical distribution in $[0,1]$, but this limiting distribution is not Turing computable. In particular, the Monte Carlo method cannot be applied to study these dynamical systems.  
\end{abstract}

\maketitle

\input{intro}
\input{preliminaries}

\input{proof}

\input{conclusion}

\newpage

\bibliographystyle{amsalpha}
\bibliography{bib_RY}

\end{document}

%% file: intro.tex

\section{Introduction}

For all practical purposes, the world around us is not a deterministic one.  Even if a simple physical system can be described deterministically,  say by the laws of Newtonian mechanics, the differential equations expressing these laws typically cannot be solved explicitly. This means that predicting the exact evolution of the system is impossible. A classical example is the famous 3-body Problem, which asks to describe the evolution of a system in which three celestial bodies (the ``Earth'', the ``Sun'', and the ``Moon'') interact with each other via the Newton's force of gravity. Computers are generally not of much help either: of course, a system of ODEs can be solved numerically, but the solution will inevitably come with an error due to round-offs. Commonly, solutions of dynamical systems are very sensitive to such small errors (the phenomenon known as ``Chaos''), so the same computation can give wildly different numerical results.

An extreme example of the above difficulties is the art of weather prediction. A realistic weather model will have such a large number of inputs and parameters that simply running a numerical computation will require a massive amount of computing resources; it is, of course, extremely sensitive to errors of computation. A classical case in point is the Lorenz system suggested by meteorologist Edward Lorenz in 1963 \cite{Lorenz}. It has only three variables and is barely non-linear (just enough not to have an explicit solution), and nevertheless it possesses a chaotic attractor \cite{Tucker} -- one of the first such examples in history of mathematics-- so deterministic weather predictions even in such a simplistic model are practically impossible.

Of course, this difficulty is well known to practitioners, and yet weather predictions are somehow made, and sometimes are even accurate. They are made in the language of statistics (e.g. there is a 40\% chance of rain tomorrow), and are based on what is broadly known as {\it Monte Carlo} technique, pioneered by Ulam and von Neumann in 1946 \cite{URvN,MetUl,Metropolis}. Informally  speaking, we can throw random darts to select a large number of initial values; run our simulation for the desired duration for each of them; then statistically average the outcomes. We then expect these averages to reflect the true statistics of our system. To set the stage more formally, let us assume that we have a discrete-time dynamical system
$$f:D\to D,\text{ where }D\text{ is a finite domain in }\mathbb R^n$$
that we would like to study. Let $\bar x_1,\ldots,\bar x_k$ be $k$ points in $D$ randomly chosen, for some $k>>1$ and consider the probability measure
\begin{equation}\label{eq:m-csum}
  \mu_{k,n}=\frac{1}{kn}\sum^k_{l=1}\sum_{m=1}^n \delta_{f^{\circ m}(\bar x_l)},
\end{equation}
where $\delta_{\bar x}$ is the delta-mass at the point $\bar x\in\mathbb R^n$. The mapping $f$ can either be given by mathematical formulas, or stand for a computer program we wrote to simulate our dynamical system. The standard postulate is then  that for $k,n\to\infty$ the probabilities $\mu_{k,n}$ converge to a limiting statistical distribution that we can use to make meaningful long-term {\it statistical} predictions of our system.

Let us say that a measure $\mu$ on $D$ is a {\it physical measure} of $f$ if its basin $B(\mu)\subset D$ --that is, the set of initial values $\bar x$ for which the weak limit of $\frac{1}{n}\sum_{m=1}^n \delta_{f^{\circ m}(\bar x)}$ equals $\mu$-- has positive Lebesgue measure. This means that the limiting statistics of such points will appear in the averages (\ref{eq:m-csum}) with a non-zero probability. If there is a unique physical measure in our dynamical system, then one random dart in (\ref{eq:m-csum}) will suffice. Of course, there are systems with many physical measures. For instance, Newhouse \cite{Newhouse} showed that a polynomial map $f$ in dimension $2$ can have infinitely many attracting basins, on each of which the dynamics will converge to a different stable periodic regime. This in itself, however, is not necessarily an obstacle to the Monte-Carlo method, and indeed, the empirical belief is that it still succeeds in these cases. 

Our results are most surprising in view of the above computational statistical paradigm. Namely we consider the simplest examples of non-linear dynamical systems: quadratic maps of the interval $[-1,1]$ of the form
$$f_a(x)=ax(1-x),\; a\in(0,4]$$
and find an uncountable set of values of $a$ for which:
\begin{enumerate}
\item there exists a {\it unique} physical measure $\mu$ and its basin $B(\mu)\subset [0,1]$ has full Lebesgue measure. 
  \item the measure $\mu$ is not computable relative to $a$. 

\end{enumerate}
This means that there is no algorithm that correctly computes $\mu$, even if the parameter $a$ is assumed to be provided to the algorithm at no computational cost.  Thus, the Monte-Carlo computational approach fails spectacularly for truly simple maps -- not because there are no physical measures, or too many of them, but because the ``nice'' unique limiting statistics cannot be computed, and thus the averages (\ref{eq:m-csum}) will not converge to anything meaningful in practice.

It is worth drawing a parallel with our recent paper \cite{YRAdvances19}, in which we studied the computational complexity of topological attractors of maps $f_a$. Such attractors  capture the limiting {\it deterministic} behavior of the orbits. They are always computable, and we found that for almost every parameter $a$, the time complexity of computing its attractor is polynomial, although there exist attractors with an arbitrarily high computational complexity. In dynamics, both in theory and in practice, it is generally assumed that long-term statistical properties are simpler to analyze than their deterministic counterparts. From the point of view of computational complexity, this appears to be false.

We note that computability of invariant measures has been studied before \cite{Cristobal, GalHoyRoj09, GaHoRo, BBRY}. In \cite{GalHoyRoj09} for instance the authors construct continuous maps of the circle for which computable invariant measures do not exists. In the context of symbolic systems, there have been some recent works studying the computational properties of the limiting statistics, see e.g. \cite{Sablik}, and of thermodynamic invariants (see e.g. in \cite{HM, BSW}). The computational complexity of individual trajectories in Hamiltonian dynamics has been addressed in e.g. \cite{kawa}. Long-term unpredictability is generally associated with dynamical systems containing embedded Turing machines (see e.g. the works \cite{moore, MK,koiran, BGR, BRS}). Dynamical properties of Turing machines viewed as dynamical systems have similarly been considered (cf. \cite{Kurka, jeandel}). Yet we are not aware of any studies of the limiting statistics in this latter context.  We also point out that a different notion of statistical intractability in dynamics, based on the complexity of a mathematical description of the set of limiting measures, has been introduced and studied in \cite{berger1,berger2}. 

From a practical point of view, some immediate questions arise. Our examples are rare in the one-parameter quadratic family $f_a(x)=ax(1-x)$. However, there are reasons to expect that in more complex multi-parametric, multi-dimensional families they can become common. Can they be {\it generic} in a natural setting? As the results of \cite{berger2} suggest, the answer may already be "yes" for quadratic polynomial maps in dimension two.
Furthermore, even in the one-dimensional quadratic family $f_a$ it is natural to ask what the typical computational complexity of the limiting statistics is -- even if it is computable in theory, it may not be in practice.  

%% file: preliminaries.tex

\section{Preliminaries}

 \subsection*{Statistical simulations and computability of probability measures}
 
We give a very brief summary of relevant notions of Computability Theory and Computable Analysis. For a more in-depth
introduction, the reader is referred to e.g. \cite{BY-book}.
As is standard in Computer Science, we formalize the notion of
an algorithm as a {\it Turing Machine} \cite{Tur}.  
We will call a  function $f:\NN\to\NN$  \emph{computable} (or {\it recursive}), if there exists a Turing Machine  $\cM$ which, upon input $n$, outputs $f(n)$.
Extending algorithmic notions to functions of real numbers was pioneered by Banach and Mazur \cite{BM,Maz}, and
is now known under the name of {\it Computable Analysis}. 
Let us begin by giving the modern definition of the notion of computable real
number,  which goes back to the seminal paper of Turing \cite{Tur}. By identifying $\Q$ with $\N$ through some effective enumeration, we can assume algorithms can operate on $\Q$. Then a real number $x\in\RR$ is called \defin{computable} if there is an algorithm  $M$ which, upon input $n$, halts and outputs a rational number $q_n$ such that  $|q_n-x|<2^{-n}$.
Algebraic numbers or  the familiar constants such as $\pi$, $e$, or the Feigenbaum constant  are computable real numbers. However, the set of all computable real numbers $\RR_C$ is necessarily countable, as there are only countably many Turing Machines. 

We now define computability of functions over $[0,1]$.  Recall that for a continuous function $f\in C_{0}([0,1])$, a modulus of continuity consists of a function $\delta: \Q\cap(0,a)\to\Q\cap(0,a)$ such that $|f(x)-f(y)|\leq \epsilon$ whenever $|x-y|\leq \delta(\epsilon)$. A function $f:[0,1]\to[0,1]$ is \defin{computable} if it has a computable modulus of continuity and there is an algorithm which, provided with a rational number which is $\delta(\epsilon)$-close to $x$, outputs a rational number which is $\epsilon$-close to $f(x)$.  

Computability of probability measures, say over $[0,1]$ for instance, is defined by requiring the ability to compute the expected value of computable functions. 

\begin{defn}\label{1}Let $(f_{i})$ be any sequence of uniformly computable functions over $[0,1]$. A probability measure $\mu$ over $[0,1]$ is \defin{computable} if there exist a Turing Machine $M$ which on input $(i,\epsilon)$ (with $\epsilon \in \Q$) outputs a rational $M(i,\epsilon)$ satisfying
$$
|M(i,\epsilon) - \int f_{i}\,d\mu | < \epsilon.
$$
\end{defn}

We note that this definition it compatible with the notion of weak convergence (see Section \ref{stage}) of measures in the sense that a measure is computable if and only if it can be algorithmically approximated (in the weak topology) to an arbitrary degree of accuracy by measures supported on finitely many rational points and with rational weights.  
Moreover, this definition also models well the intuitive notion of statistical sampling in the sense that a measure $\mu$ is computable if and only if there is an algorithm to convert sequences sampled from the uniform distribution into sequences sampled with respect to $\mu$. 

In this paper, we will be interested in the computability properties of invariant measures of  quadratic maps of the form $ax(1-x)$, with $a\in \RR$. As is standard in computing practice, we will assume that the algorithm can read the value of $a$ externally in order to compute $\mu$. More formally, let us denote $\cD_n\subset \RR$ the set of dyadic rational numbers with denominator $2^{-n}$. 
We say that a function $\phi:\NN\to\QQ$ is an {\it oracle} for $a\in \RR$ if for every $m\in \NN$
$$\phi(m)\in \cD_m\text{ and }d(\phi(m),a)<2^{-(m-1)}.$$
We amend our definitions of computability of a probability measure $\mu$ by allowing {\it oracle Turing Machines} $M^\phi$ where $\phi$ is any function as above. 
On each step of the algorithm, $M^\phi$ may read the value of $\phi(m)$ for an arbitrary $m\in\NN$.
This approach, usually referred to as computability {\it relative to} $a$, allows us to separate the questions of computability of a parameter $a$ from that of the measure.


\subsection*{Invariant measures of quadratic polynomials and the statement of the main result.}
As before, we denote 
$$f_a(x)=ax(1-x).$$
For $a\in[0,4]$, this quadratic polynomial maps the interval $[0,1]$ to itself. We will view $f_a:[0,1]\to[0,1]$ as a discrete dynamical system, and will denote $f_a^n$ the $n$-th iterate of $f_a$. 

A measure $\mu$ is called \defin{physical} or {\it Sinai-Ruelle-Bowen (SRB)} if 
\begin{equation}\label{srb}
\frac{1}{n}\sum_{k=0}^{n-1}\delta_{f^{k}x} \to \mu 
\end{equation}
for a set of positive Lebesgue measure.  It is known that if a physical measure exists for a quadratic map $f_{a},\,\, a\in [0,4]$, then it is unique and (\ref{srb}) is satisfied for Lebesgue almost all $x\in [0,1]$.

\bigskip
\noindent\textbf{Main Theorem.} \emph{
There exists parameters $a\in(0,4)$ for which the quadratic map $f_a(x)=ax(1-x)$ has a physical measure $\mu$ which is not computable relative to $a$. 
	}
\bigskip

%% file: proof.tex

\section{Proof of the Main Theorem}
The proof is based on a delicate construction in one-dimensional dynamics described in \cite{HofKel}, which will allow us to construct maps $f_a$ with physical measures which selectively charge points in a countable set of periodic orbits. To give precise formulations, we will need to introduce some further concepts.

\subsection{Setting the stage}\label{stage}
It will be convenient to recall that weak convergence of measures on $[0,1]$ is compatible with the notion of {\it Wasserstein-Kantorovich distance},
defined by:

\begin{equation*}
W_{1}(\mu,\nu)=\underset{f\in 1\text{-Lip}([0,1])}{\sup}\left|\int f d\mu-\int f d\nu\right|
\end{equation*}
\noindent where $1\mbox{-Lip}([0,1])$ is the space of Lipschitz functions on $[0,1]$, having Lipschitz constant less than one.

For $a\in[0,4]$ and $x\in[0,1]$, we set
$$\nu_a^n(x)=\frac{1}{n}\sum_{k=0}^{n-1}\delta_{f_a^kx}.$$

We will make use of the following folklor fact (see e.g. \cite{dMvS}):
\begin{prop}
  \label{sink}
  Suppose, for $a\in[0,4]$ the map $f_a$ has an attracting periodic orbit of period $p$:
  $$x_0\overset{f_a}{\mapsto}x_1\overset{f_a}{\mapsto}\cdots\overset{f_a}{\mapsto}x_{p-1}\overset{f_a}{\mapsto}x_0,\;\left| \frac{d}{dx}f_a^p(x_0)\right|<1.$$
  Let
  $$\mu\equiv \frac{1}{p}\sum_{k=0}^{p-1}\delta_{x_k}.$$
  Then $\mu$ is the unique physical measure of $f_a$ (so, in particular, the attracting orbit is unique); and
  $$W_1(\nu_a^n(x), \mu)\to 0$$
  uniformly on a set of full Lebesgue measure in $[0,1]$.

  \end{prop}

For $a\in(0,4]$ consider the third iterate $g_a\equiv f_a^3$. 
We start by noting that there exists a parameter value $c\in (3.85,4)$ such that the following holds: $$g_c(0.5)\neq g_c^2(0.5)=g_c^3(0.5).$$
If we denote $\beta_c=g_c^2(0.5)$, then $\beta'_c\equiv g_c(0.5)=1-\beta_c$, and denoting $I_c\equiv [\beta'_c,\beta_c]\ni 0.5$, we have
$$g_c(I_c)=I_c;$$
both endpoints of $I_c$ map to $\beta_c$. The restriction $g_c|_{I_c}$ maps both halfs ($L_c=[\beta',0.5]$ and $R_c=[0.5,\beta]$)  of the  interval $I_c$ onto the whole $I_c$ in a monotone fashion (that is, it folds $I_c$ over itself).

\begin{figure}[ht]
  \includegraphics[width=0.6\textwidth]{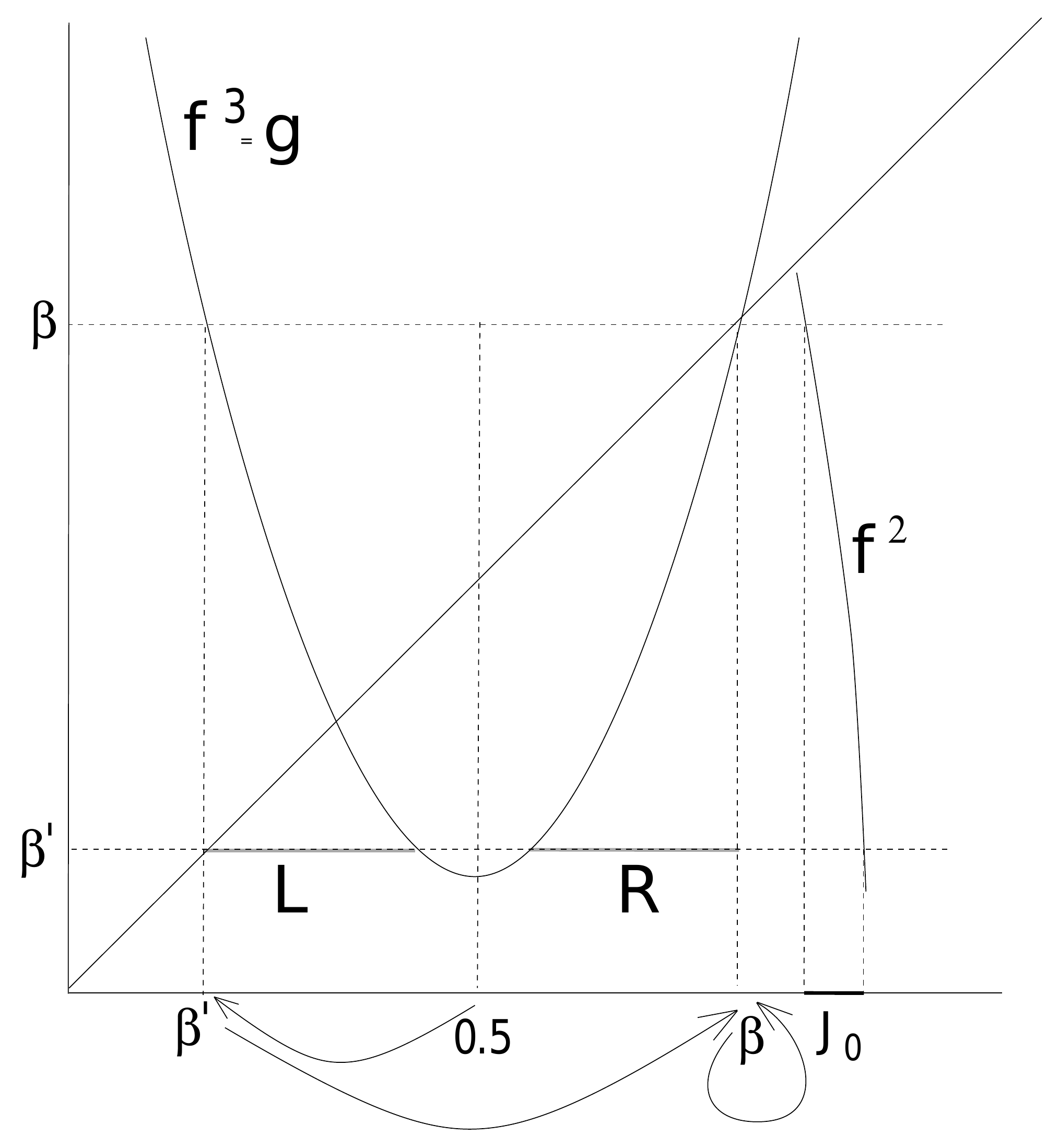}
\caption{\label{fig:iterates}Some iterates of $f\equiv f_a$ for $a\in(c,4]$ (we drop the subscript $a$ for simplicity in all notations in the figure).}
 \end{figure} 

For $a\in [c,4]$, there exists a continuous branch $\beta_a$ of the fixed point
$$g_a(\beta_a)=\beta_a,$$
and we again set $\beta'_a=1-\beta_a$  (so $g_a(\beta'_a)=\beta_a$), and $I_a\equiv[\beta'_a,\beta_a]$. Now, if $a\in(c,4]$, the image
  $$g_a(I_a)\supsetneq I_c,\text{ with }g_a(0,5)<\beta'_a.$$
  Thus, there is a pair of sub-intervals $L_a=[\beta'_a,l_a]$, $R_a=[r_a,\beta_a]$ inside $I_a$ which are mapped monotonely over $I_a$ by $g_a$ (the endpoints $l_a$, $r_a$ are both mapped to $\beta'_a$. See Figure~\ref{fig:iterates} for an illustration.

  Assigning values $0$ to $L_a$, and $1$ to $R_a$ we obtain symbolic dynamics on the set of points
  $$C_a\equiv \{x\in I_a\text{ such that }g_a^n(x)\in I_a\text{ for all }n\in\NN\}.$$
  If $a=c$, then, of course, $C_a=I_a$. Otherwise, the following is well-known:
  \begin{prop}\label{prop-cantor}
If $a\in(c,4)$ then $C_a$ is a Cantor set, and the symbolic dynamics conjugates $g_a|_{C_a}$ to the full shift on $\{0,1\}$.
    \end{prop}
  In particular, every periodic sequence of $0$'s and $1$'s corresponds to a unique periodic orbit in $C_a$ with this symbolic dynamics. These orbits clearly move continuously with $a$, and can be easily computed given $a$ and the symbolic sequence as the unique fixed points of the corresponding monotone branches of iterates of $g_a$.

  We  enumerate all periodic sequences of $0's$ and $1$'s as follows. A sequence with a smaller period will precede a sequence with a larger period. Within the sequences of the same period, the ordering will be lexicographic, based on  
  the convention $1\prec 0$. 
We let $$\{p^1_{a,n},\ldots,p^{k_n}_{a,n}\}$$ be the periodic orbit of $g_a$ in $C_a$ which corresponds to the $n$-th symbolic sequence in this ordering (note that the first one is $\beta_a$). We denote
  $$\text{Per}_a(n)=\cup_{j=0}^2 f^j_a(\{p^1_{a,n},\ldots,p^{k_n}_{a,n}\}),$$
  which is, clearly, a periodic orbit of $f_a$.
Let us denote $$\lambda_a(n)=\frac{1}{|\text{Per}_a(n)|}\sum_{x\in\text{Per}_a(n)}\delta_x.$$

  \subsection{Main construction}
  Our arguments will be based on the results of  F.~Hofbauer and G.~Keller in \cite{HofKel}; see also the earlier paper of S.~Johnson \cite{Johnson}, which uses similar language to ours.

  Let us develop some further notation. For $x\in [0,1]$ and $a\in[0,4]$ we let
  $\Omega_a(x)$ denote the set of weak limits of the sequence $\nu_a^n(x)=\frac{1}{n}\sum_{k=0}^{n-1}\delta_{f_a^kx}$. Let us denote $\cP\subset(c,4]$ the collection of parameters $a$ such that the following holds:
    \begin{itemize}
    \item $f_a$ has a unique physical probability measure $\mu_a$;
    \item denoting $m_a(n)\equiv\mu_a(\text{Per}_a(n))$, we have
      \begin{equation}  \label{eq-charge}
\sum_{n=1}^\infty m_a(n)=1.
      \end{equation}

      \end{itemize}
    Thus, the charge of the physical measure $\mu_a$ resides in the periodic orbits in the Cantor set $C_a$.

    We will formulate the following direct consequence of the main result (Theorem~5) of \cite{HofKel}\footnote{Note that the set of physical measures constructed in Theorem 5 of \cite{HofKel} includes convex combinations of $\lambda_a(n)$. Compare also with Theorem~1 of \cite{HofKel}}:
    \begin{thm}
      \label{th:hk}
      There exists an infinite set $\hat \cP\subset\cP$ such that the following holds:
      \begin{enumerate}
      \item $\Omega_a(x)=\Omega_a(0.5)=\{\mu_a\}$ for Lebesgue almost every $x$;
      \item for any sequence of non-negative reals $\{l_n\}_{n\in\NN}$ with $\sum l_n=1$, the subset of $a\in\hat \cP$ for which
        $m_a(n)=l_n$ is dense in $\hat \cP$.
        \end{enumerate}
      \end{thm}

    We will outline the idea of the construction of such maps below, but the complete proof of Theorem~\ref{th:hk} is quite technical and goes beyond the scope of this paper.


    We start with the following ``simple'' example:

    \medskip\noindent{\bf Example 1}: {\it The set $$\Omega_a(0.5)=\lambda_a(1)=\frac{1}{3}(\delta_{\beta_a}+\delta_{f_a(\beta_a)}+\delta_{f^2_a(\beta_a)})$$ and
      $\Omega(x)=\Omega(0.5)$ for almost every $x$ (compare with Theorem~1 of \cite{HofKel}).}

    \medskip\noindent
    Consider again Figure~\ref{fig:iterates} as an illustration. We note that there exists an interval $J_0$ to the right of the fixed point $\beta_a$ such that the following holds:
    \begin{itemize}
    \item $f^2_a(J_0)\Supset [\beta_a',\beta_a]$;
      \item Denote by $\psi_a$ the branch of $g_a^{-1}$ which fixes $\beta_a$. Then the interval $J_0$ is contained in the domain of definition of $\psi_a$. Thus, there is an orbit $$J_{-n}\equiv \psi_a^n(J_0) \to\beta_a \;(\text{here }f_a^{3n}(J_{-n})=J_0).$$
      \end{itemize}
    Moving the parameter $a\in(c,4]$, we can place the image $g_a^2(0.5)$ at any point of  $J_{-{n_1}}$, for an arbitrary $n_1$. If the value of $n_1$ is large, then the $g_a$-orbit of $0.5$ will spend a long time in a small neighborhood of  $\beta_a$, before hitting some $x_1\in J_0$.
  Adjusting the position of
      $a\in(c,4]$, we can ensure that $f_a^2(x_1)$ is inside $J_{-n_2}$ for an even larger $n_2$, so the orbit returns to an even smaller neighborhood of $\beta_a$ where it will spend an even longer time. Continuing increasing $n_k$'s as needed so the orbit of $0.5$ spends most of its time in ever smaller neighborhoods of $\beta_a$, we can ensure that the averages
    $\nu_a^n(0.5)=\frac{1}{n}\sum_{k=0}^{n-1}\delta_{f_a^k(0.5)} $ converge to the delta masses supported on the orbit of $\beta_a$.

   Proceeding in this way,
    for an arbitrarily large $l\in\NN$ and $x\in [\beta_a',\beta_a]$ we can find $a\in(c,4]$ and $m>2^{-l}$ such that:
    \begin{enumerate}
    \item the distance
      $$W_1(\nu_a^{m}(0.5)-\lambda_a(1))<2^{-l};$$
    \item the iterate $f_a^m(0.5)$ lies in $J_0$;
    \item the next iterate $f_a^{m+1}(0.5)=0.5$.  

    \end{enumerate}

    Property (3) ensures that the critical point $0.5$ is periodic with period $m+1$. Since $f_a'(0.5)=0$, we have
		$(f_a^{m+1})'(0.5)=0$, so this is a (super)attracting periodic point. Proposition~\ref{sink} implies that the physical measure $\mu_a$ for $f_a$ is supported on the orbit of $0.5$, and thus
    $$W_1(\mu_a-\lambda_a(1))<2\cdot 2^{-l}.$$
    Again, by Proposition~\ref{sink} and considerations of continuity, there exist $n>m$ and $\eps>0$ such that for any $a'$ with $|a'-a|<\eps$, we have
    $$W_1(\nu_{a'}^{n}(x)-\lambda_a(1))<4\cdot 2^{-l}$$
    for any $x$ in a set of length $1-2^{-l}$.

    Assuming $\eps$ is small enough, we again have $g_{a'}(0.5)$ slightly to the right of $\beta_{a'}$ and we can repeat the above steps inductively to complete the construction.
    
        As a next step, we  construct an asymptotic  measure supported on two periodic orbits:

        \medskip
        \noindent
            {\bf Example 2:} {\it the set $\Omega(0.5)=a_1\lambda_1(1)+a_2\lambda_a(n)$ for $n>1$ and $a_1+a_2=1$.}

            \medskip\noindent
            Let $p\in \text{Per}_a(n)$ and, as before, denote by $3k_n$ its period. Letting $\phi$ denote the branch of $f_a^{-3k_n}$ fixing $p$,
            we again find a $\phi$-orbit
            $$J'_0\mapsto J'_{-1}\mapsto J'_{-2}\mapsto\cdots,\text{ with }J'_{-k}\to p,$$
            where
            $$f_a^s(J_0)\Supset [\beta'_a,\beta_a]$$
            for a univalent branch of the iterate $f_a^s$.

            Now we can play the same game as in Example 1, alternating between entering the orbit $J_{-k}$ close to the point $\beta_a$, and the orbit $J'_{-k}$ close to $p$. In this way, we can achieve the desired limiting asymptotics with any values $a_1, a_2$.

The above construction can be clearly modified for any countable collection of periodic orbits in $C_a$, as required for the proof of Theorem~\ref{th:hk}.

\subsection{Constructing non computable physical measures}
\begin{definition}
  Let us define a very particular subset $\tl \cP\subset \hat\cP$ as follows: $a\in\tl\cP$ if
\begin{equation}
m_a(2n-1)+m_a(2n)=2^{-n}   \qquad \text{ for all } n\in \N.
\end{equation}
\end{definition}

For convenience of reference, let us formulate a corollary of Theorem~\ref{th:hk}:

\begin{prop}\label{induction.step} Suppose, $a\in\tl\cP$. 
  Then, for every $\epsilon>0$, $l\in\N$ and $s\in\{0,1\}$, there exists $a'>a$ such that
  \begin{itemize}
\item $a'\in\tl\cP$;   
\item $|a-a'|<\epsilon$;
\item $m_a(n)=m_{a'}(n)$ for all $n\notin \{2l, 2l-1\}$ and
\item $m_{a'}(2l-s)=2^{-l}$. 
\end{itemize}
\end{prop}

Let $(\tau_{i})_{i\in\NN}$ be the smallest collection of functions containing the \emph{step} continuous functions of rational intervals, and which is closed by rational linear combinations and scalar multiplication. Note that this is a countable collection of functions that can be enumerated in an effective way.  

We construct a parameter $a$ for which the map $f_a=ax(1-x)$ has a unique physical measure $\mu_a$ such that for any Turing Machine $M^{\phi}$ with an oracle $\phi$ for $a$, that computes a probability measure, there exists $i$ and $\epsilon>0$ such that 

$$
|M^{\phi}(i,\epsilon) - \int \tau_{i}\, d\mu|  > \epsilon.
$$

Our construction can be thought of as a game between a \emph{Player} and infinitely many \emph{opponents}, which will correspond to the sequence consisting of
machines $M^{\phi}_n$ that compute some probability measure. The opponents try to compute $\mu_a$ by asking the Player to provide an oracle $\phi$ for $a$,
while the Player tries to chose the bits of $a$ in such a way that none of the opponents correctly computes $\mu_a$.

We show that the Player always has a winning strategy: it plays against each machine, one by one, asking the machine to compute the expected value of a particular function $\tau_{i}$ to a certain degree of accuracy. The machine then runs for a while, asking the Player to provide more and more bits of $a$, until it eventually halts and outputs a rational number. Then the Player reveals the next bit of $a$ and shows that the machine's answer is incompatible with $\mu_{a}$. The details are as follows.


We will proceed inductively.  Let $M^{\phi}_{1}, M^{\phi}_{2},\dots$ be some enumeration of all the machines with an oracle for $a$ that compute some probability measure. 
 At step $n$ of the induction, we will have a parameter $a_n\in(c,4)$ and a natural number $l_{n}$ such that:
\begin{enumerate}
\item $a_n\in\tl\cP$;
\item\label{fooling} there exists $i=i(n)\in \N$ such that either
\begin{itemize}
\item $M_{n}^{\phi}(\tau_{i},2^{-{n}}/100) \leq 2^{-{n}}/2$  whereas $\mu_{a_{n}}(\tau_{i})\sim2^{-{n}}$; or
\item $M_{n}^{\phi}(\tau_{i},2^{-{n}}/100) > 2^{-{n}}/2$  whereas $\mu_{a_{n}}(\tau_{i})\sim0$
\end{itemize}
In other words, given an oracle for $a_n$, the machine $M_n^\phi$ cannot correctly approximate the value of $\mu_{a_{n}}$ at $\tau_{i}$;
\item\label{measures}  $|\mu_{a_{n-1}}(\tau_{i(k)}) - \mu_{a_{n}}(\tau_{i(k)})| < 2^{-3n}$ for all $k<n$; 
\item\label{parameters} $|a_n-a_{n-1}|<2^{-3l_n}$.
\end{enumerate}
\bigskip

\noindent\textbf{Base of the induction.}  We start by letting $a$ be any of the parameters in $\tl \cP$. We note that 
$
m_a(1)+m_a(2)=2^{-1}.
$
 It follows that there exists $\tau=\tau_{i(1)}$ such that $
|\mu_{a}(\tau) - m_a(1)|<2^{-1}/200$\footnote{Note that the mass of the higher periodic points that may fall in an open set containing $\text{Per}_a(1)$ goes to zero as the diameter of the open set goes to zero.}, and 

\begin{equation}
\label{tau} \supp \tau \cap \text{Per}(j) = \emptyset \,\, \text{ for all } \,\, 1<j<10.  
\end{equation}

We now let the machine $M^{\phi}_{1}$ compute the expected value of $\tau$ with precision $2^{-1}/100$, giving it $a$ as the parameter. Let $l_1$ be the last time a bit of $a$ is queried by $M^{\phi}_{1}$ during the computation. By Proposition \ref{induction.step}, for any $s\in \{0,1\}$ there exists $a'$ such that 
\begin{itemize}
\item $|a-a'|<2^{-3l_1}$;
\item $a'\in\tl\cP$;
\item $m_{a'}(2-s)=2^{-1}$. 
\end{itemize}

Let $q=M^{\phi}_{1}(\tau,2^{-1}/100)$.  There are two possibilities:

\begin{enumerate}
\item[Case 1. ]\label{case1} If $q\leq 2^{-1}/2$, we chose $a'$ above so as to have $m_{a'}(1)=2^{-1}$;
\item[Case 2. ]\label{case2} If $q>2^{-1}/2$, we chose $a'$ above so as to have $m_{a'}(2)=2^{-1}$ (and therefore $m_{a'}(1)=0$);
\end{enumerate}
We then let $a_{1}\equiv a'$. By \ref{tau}, $|\mu_{a_1}(\tau)-m_{a_1}(1)|<2^{-10}$. Note that up to the first $l$ bits, $a_{0}$ and $a_{1}$ are indistinguishable and therefore the machine $M^{\phi}_{1}$ will return the same answer for both parameters. It follows that the machine $M^{\phi}_{1}$ cannot correctly approximate $\mu_{a_1}$ at $\tau$. 

\bigskip

\noindent\textbf{Step of the induction.} Assume $a_{n-1}\in \tl \cP$ has been defined. Then it holds 
$$
m_{a_{n-1}}(2n-1)+m_{a_{n-1}}(2n) = 2^{-n},
$$
and there exists $\tau=\tau_{i(n)}$ such that $|\mu_{a_{n-1}}(\tau) - m_{a_{n-1}}(2n-1)|<2^{-n}/200$ and
\begin{equation}
\label{taun} \supp \tau \cap \text{Per}(j) = \emptyset \,\, \text{  for all }\,\, 2n-1<j<10n.   
\end{equation}

Once again, we let the machine $M^{\phi}_{n}$ compute the expected value of $\tau$ with precision $2^{-n}/100$, giving it $a_{n-1}$ as the parameter. Let $l_n$ be the last time a bit of $a_{n-1}$ is queried by $M^{\phi}_{1}$ during the computation. By Proposition \ref{induction.step} again, for any $s\in \{0,1\}$ there exists $a'$ such that
\begin{itemize}
\item $|a_{n-1}-a'|<2^{-3l_n}$;
\item $a'\in\tl\cP$;
\item $m_{a_{n-1}}(t)=m_{a'}(t)$ for all $t\notin \{2n-1,2n\}$ and
\item $m_{a'}(2n-s)=2^{-n}$. 
\end{itemize}

Let $q=M^{\phi}_{n}(\tau,2^{-n}/100)$.  There are two possibilities:

\begin{enumerate}
\item[Case 1. ]\label{case1} If $q\leq 2^{-n}/2$, we chose $a'$ above so as to have $m_{a'}(2n-1)=2^{-n}$;
\item[Case 2. ]\label{case2} If $q>2^{-n}/2$, we chose $a'$ above so as to have $m_{a'}(2n)=2^{-n}$ (and therefore $m_{a'}(2n-1)=0$);
\end{enumerate}
We then let $a_{n}\equiv a'$. Since $\tau$ satisfies property \ref{taun}, we have that $|\mu_{a_n}(\tau)-m_{a_n}(2n-1)|<2^{-10n}$.  Note that up to the first $l_n$ bits, $a_{n-1}$ and $a_{n}$ are indistinguishable, and thus the machine $M^{\phi}_{n}$ will return the same answer for both parameters. It follows that the machine $M^{\phi}_{n}$ cannot correctly approximate $\mu_{a_n}$ at $\tau$.  Moreover, property \ref{taun} again and the fact that (by construction) $a_n$ satisfies   
$$
m_{a_{n-1}}(t)=m_{a_n}(t) \,\,\text{ for all } \,\, t\notin \{2n-1,2n\},
$$
guarantee that Condition (\ref{measures}) is satisfied as well. We now let $a_{\infty}=\lim_na_n$ and claim that $\mu_{a_{\infty}}$ has the required properties. Indeed, Condition (\ref{fooling}) ensures that for every $n$ there is a step function $\tau_{i(n)}$ at which machine $M^{\phi}_n$ fails to compute correctly the expected value for $\mu_{a_n}$, and Condition (\ref{measures}) guarantees that the same holds for $\mu_{a_{\infty}}$.

%% file: conclusion.tex
\section{Conclusion}
\label{sec:conclusion}

Ever since the first numerical studies of chaotic dynamics appeared in the early 1960's (such as the work of Lorenz \cite{Lorenz}), it has become commonly accepted among practitioners that computers cannot, in general, be used to make deterministic predictions about future behavior of nonlinear dynamical systems. Instead, the standard practice now is to make statistical predictions. This approach is based on the Monte Carlo method, pioneered by Ulam and von Neumann at the dawn of the computing age. It is universal and powerful -- and only requires access to the dynamical system as a black box, which is then subjected to a number of statistical trials.
Applications of the Monte Carlo technique are ubiquitous, ranging from weather forecasts to simulating nuclear weapons tests (nuclear weapons design was, of course, the original motivation of its inventors).

Our result raises a disturbing possibility that even for the most simple family of examples of non-linear dynamical systems the Monte Carlo method can fail. Given one of our examples as a black box, no algorithm can find its limiting statistics. How common such examples are in higher-dimensional families of dynamical systems, and whether one is likely to encounter one {\it in practice} remain exciting open questions.

%% file: main.bbl
\providecommand{\bysame}{\leavevmode\hbox to3em{\hrulefill}\thinspace}
\providecommand{\MR}{\relax\ifhmode\unskip\space\fi MR }
\providecommand{\MRhref}[2]{%
  \href{http://www.ams.org/mathscinet-getitem?mr=#1}{#2}
}
\providecommand{\href}[2]{#2}
\begin{thebibliography}{BBRY11}

\bibitem[BB19]{berger2}
P.~Berger and J.~Bochi, \emph{On emergence and complexity of ergodic
  decompositions}, ar{X}iv preprint math/1901.03300 (2019).

\bibitem[BBRY11]{BBRY}
I.~Binder, M.~Braverman, C.~Rojas, and M.~Yampolsky, \emph{Computability of
  {B}rolin-{L}yubich measure}, Commun. Math. Phys. \textbf{308} (2011),
  743--771.

\bibitem[Ber17]{berger1}
P.~Berger, \emph{Unpredictability of dynamical systems and non-typicality of
  the finiteness of the number of attractors in various topologies}, Tr. Mat.
  Inst. Steklova \textbf{297} (2017), 7--37, English version published in Proc.
  Steklov Inst. Math. {{\bf{297}}} (2017), no. 1, 1--27.

\bibitem[BGR12]{BGR}
Mark Braverman, Alexander Grigo, and Cristobal Rojas, \emph{Noise vs
  computational intractability in dynamics}, Proceedings of the 3rd Innovations
  in Theoretical Computer Science Conference (New York, NY, USA), ITCS '12,
  ACM, 2012, pp.~128--141.

\bibitem[BM37]{BM}
S.~Banach and S.~Mazur, \emph{Sur les fonctions caluclables}, Ann. Polon. Math.
  \textbf{16} (1937).

\bibitem[BRS15]{BRS}
Mark Braverman, Cristobal Rojas, and Jon Schneider, \emph{Space-bounded
  {C}hurch-{T}uring thesis and computational tractability of closed systems},
  Physical Review Letters \textbf{115} (2015), no.~9.

\bibitem[BW18]{BSW}
M.~Burr, M.;~Schmoll and C.~Wolf, \emph{On the computability of rotation sets
  and their entropies.}, Ergodic Theory and Dynamical Systems (2018), 1--35.

\bibitem[BY08]{BY-book}
M~Braverman and M.~Yampolsky, \emph{Computability of {J}ulia sets}, Algorithms
  and {C}omputation in {M}athematics, vol.~23, Springer, 2008.

\bibitem[dMvS93]{dMvS}
W.~de~{M}elo and S.~van {Strien}, \emph{One-dimensional dynamics},
  Springer-{V}erlag, 1993.

\bibitem[GHR10]{GalHoyRoj09}
S.~Galatolo, M.~Hoyrup, and C.~Rojas, \emph{Dynamics and abstract
  computability: computing invariant measures}, Discr. Cont. Dyn. Sys. Ser A
  (2010).

\bibitem[GR11]{GaHoRo}
M.~G{\'a}cs, P.;~Hoyrup and C.~Rojas, \emph{Randomness on computable
  probability spaces -- a dynamical point of view}, Theory Comput. Syst.
  \textbf{48} (2011), no.~465.

\bibitem[HdMS16]{Sablik}
Benjamin Hellouin~de Menibus and Mathieu Sablik, \emph{Characterisation of sets
  of limit measures after iteration of a cellular automaton on an initial
  measure}, Ergodic Theory and Dynamical Systems \textbf{38} (2016), no.~2,
  601--650.

\bibitem[HK90]{HofKel}
F.~Hofbauer and G.~Keller, \emph{Quadratic maps without asymptotic measure},
  Comm. {M}ath. {P}hys. \textbf{127} (1990), 319--337.

\bibitem[HM10]{HM}
M.~Hochman and T.~Meyerovitch, \emph{Characterization of the entropies of
  multidimensional shifts of finite type.}, Annals of Mathematics \textbf{171}
  (2010), no.~3, 2011--2038.

\bibitem[Jea14]{jeandel}
E.~Jeandel, \emph{Computability of the entropy of one-tape {T}uring machines},
  31st {I}nternational {S}ymposium on {T}heoretical {A}spects of {C}omputer
  {S}cience (STACS), LIPIcs. Leibniz Int. Proc. Inform., vol.~25, 2014,
  pp.~421--432.

\bibitem[Joh87]{Johnson}
S.~Johnson, \emph{Singular measures without restrictive intervals}, Commun.
  {M}ath. {P}hys. \textbf{110} (1987), 185--190.

\bibitem[KCG94]{koiran}
P.~Koiran, M.~Cosnard, and M.~Garzon, \emph{Computability with low-dimensional
  dynamical systems}, Theoret. Comput. Sci. \textbf{132} (1994), no.~1-2,
  113--128.

\bibitem[KTZ18]{kawa}
Akitoshi Kawamura, Holger Thies, and Martin Ziegler, \emph{Average-case
  polynomial-time computability of {H}amiltonian dynamics}, 43rd International
  Symposium on Mathematical Foundations of Computer Science, {MFCS} 2018,
  August 27-31, 2018, Liverpool, {UK}, 2018, pp.~30:1--30:17.

\bibitem[Kur97]{Kurka}
Petr Kurka, \emph{On topological dynamics of {T}uring machines.}, Theoret.
  Comput. Sci. \textbf{174} (1997), 203--2016.

\bibitem[Lor63]{Lorenz}
E.~N. Lorenz, \emph{Deterministic nonperiodic flow}, J. Atmos. Sci. \textbf{20}
  (1963), 130--141.

\bibitem[Maz63]{Maz}
S.~Mazur, \emph{Computable {A}nalysis}, vol.~33, Rosprawy Matematyczne, Warsaw,
  1963.

\bibitem[Met87]{Metropolis}
N.~Metropolis, \emph{The beginning of the {M}onte {C}arlo method}, Los {A}lamos
  {S}cience {S}pecial {I}ssue (1987), 125--130.

\bibitem[MK99]{MK}
C.~Moore and P.~Koiran, \emph{Closed-form analytic maps in one and two
  dimensions can simulate universal {T}uring machines.}, Theoret. Comput. Sci.
  \textbf{210} (1999), no.~1, 2217--223.

\bibitem[Moo91]{moore}
C.~Moore, \emph{Generalized shifts: unpredictability and undecidability in
  dynamical systems}, Nonlinearity \textbf{4} (1991), no.~2, 199--230.

\bibitem[MS49]{MetUl}
N.~Metropolis and Ulam. S., \emph{The {M}onte {C}arlo method}, Journal of the
  {A}merican {S}tatistical {A}ssociation \textbf{44} (1949), 335--341.

\bibitem[New74]{Newhouse}
S.~Newhouse, \emph{Diffeomorphisms with infinitely many sinks}, Topology
  \textbf{13} (1974), 9--18.

\bibitem[Roj08]{Cristobal}
C.~Rojas, \emph{Randomness and ergodic theory: an algorithmic point of view},
  Ph.D. thesis, Ecole Polytechnique, 2008.

\bibitem[RY19]{YRAdvances19}
Cristobal Rojas and Michael Yampolsky, \emph{Computational intractability of
  attractors in the real quadratic family}, Advances in Mathematics
  \textbf{349} (2019), 941 -- 958.

\bibitem[Tuc02]{Tucker}
W.~Tucker, \emph{A rigorous {O}{D}{E} solver and {S}male's 14th problem},
  Found. Comp. Math. \textbf{2} (2002), 53--117.

\bibitem[Tur36]{Tur}
A.~M. Turing, \emph{On computable numbers, with an application to the
  {E}ntscheidungsproblem}, Proceedings, London Mathematical Society (1936),
  230--265.

\bibitem[URvN47]{URvN}
S.~Ulam, R.D. Richtmyer, and J.~von Neumann, \emph{Statistical methods in
  neutron diffusion}, Los {A}lamos {S}cientific {L}aboratory report
  {L}{A}{M}{S} 551 (1947).

\end{thebibliography}
